\documentclass[12pt,leqno,a4paper]{article}
\usepackage{amsmath,amsthm,calc}
\usepackage[latin1]{inputenc}
\usepackage[T1]{fontenc}
\usepackage[english]{babel}

\newtheorem{theorem}{Theorem}[section]

\newtheorem{lemma}[theorem]{Lemma}
\newtheorem{corollary}[theorem]{Corollary}
\newtheorem{definition}[theorem]{Definition}

\newtheorem{example}[theorem]{Example}
\newtheorem{remark}[theorem]{Remark}


\newcommand*\proofnamestyle{\itshape}

\DeclareMathOperator{\tr}{Tr}

\begin{document}

    \title{Monotone trace functions of several variables}
    \author{Frank Hansen\thanks{The author would like to dedicate this paper
    to the memory of Gert K. Pedersen.}}
    \date{September 27 2004}
    \maketitle

    \begin{abstract} We investigate monotone operator functions of several variables under a trace or a 
    trace-like functional. In particular, we prove
    the inequality $ \tau(x_1\cdots x_n)\le\tau(y_1\cdots y_n) $ for a trace $ \tau $
    on a $ C^* $-algebra
    and abelian $ n $-tuples $ (x_1,\dots,x_n)\le (y_1,\dots,y_n) $ of positive elements. We formulate and
    prove Jensen's inequality for expectation values, and we study matrix functions 
    of several variables which are convex or monotone with respect to the weak majorization for matrices.
    \end{abstract}

    \section{Preliminaries} The question of monotonicity of operator functions under a trace or a trace-like
    functional was considered in \cite{kn:hansen:2003:3}. In particular, conditions
    were given for which the implication
    \begin{gather}\label{monotone trace function}
    \underline{x}\le\underline{y}\quad\Rightarrow\quad\varphi(f(\underline{x}))\le \varphi(f(\underline{y}))
    \end{gather}
    is valid for a positive functional $ \varphi $ on a $ C^* $-algebra $ {\mathcal A} $ and abelian $ n $-tuples
    $ \underline{x} $ and $ \underline{y} $ in
    $ {\mathcal A} $ contained in the domain of a function $ f $ of $ n $ variables. 
    
    Consider $ n $-tuples $ \underline{x}=(x_1,\dots,x_n) $ and $ \underline{y}=(y_1,\dots,y_n) $ of elements
    in a $ C^* $-algebra $ {\mathcal A} $ and recall that $ \underline{x} $ is said to be abelian if the elements 
    $ x_1,\dots,x_n $ are mutually commuting; and we write $ \underline{x}\le\underline{y} $ if $ x_i\le y_i $
    for $ i=1,\dots,n. $ We say that $ \underline{x} $ is in the domain of a real continuous function $ f $ of
     $ n $ variables defined on a cube 
    $ \underline{I}=I_1\times\cdots\times I_n $ where each $ I_i $ is an interval and $ x_i $ is self-adjoint, if
    the spectrum $ \sigma(x_i) $ of $ x_i $ is contained in $ I_i $ for $ i=1,\dots,n. $ 
    In this situation $ f(\underline{x}) $ is
    naturally defined as an element in $ {\mathcal A}. $ 
    
    It is proved in \cite{kn:hansen:2003:3} that (\ref{monotone trace function}) holds, if $ f $ is continuous,
    convex and
    separately increasing, and the elements $ x_1,\dots,x_n $ are contained in the centralizer
    \[
    {\mathcal A}^{\varphi}=\{y\in{\mathcal A}\mid \varphi(xy)=\varphi(yx)\,\forall x\in{\mathcal A}\}.
    \] 
    
    Likewise (\ref{monotone trace function}) holds, if $ f $ is continuous, concave and separately increasing, and the 
    elements $ y_1,\dots,y_n $ are contained in the centralizer $ {\mathcal A}^\varphi. $ Finally
    (\ref{monotone trace function}) holds, if $ f $ is just continuous and separately increasing, the
    elements $ x_1,\dots,x_n $ and $ y_1,\dots,y_n $ are contained in the centralizer $ {\mathcal A}^\varphi, $
    and $ \underline{x} $ and $ \underline{y} $ are compatible in the sense that
    the commutators satisfy
    \begin{gather}\label{compatible pairs}
    \,[x_i, y_j]=[x_j, y_i]\qquad i,j=1,\dots,n.
    \end{gather}
    However, this last condition is very restrictive. It is equivalent to the demand that the line through
    $ \underline{x} $ and $ \underline{y} $ (or indeed just the midpoint) consists of abelian $ n $-tuples.
    
    It is the aim of the present article to remove the
    compatibility condition (\ref{compatible pairs}) and prove that (\ref{monotone trace function}) holds 
    for a large class of
    separately increasing functions $ f $ and elements  $ x_1,\dots,x_n $ and $ y_1,\dots,y_n $ in the centralizer 
    $ {\mathcal A}^\varphi. $ In particular, we prove the implication 
    \[
    (x_1,\dots,x_n)\le(y_1,\dots,y_n)\quad\Rightarrow\quad\tau(x_1\cdots x_n)\le\tau(y_1\cdots y_n)
    \]
    for a trace $ \tau $ on a $ C^* $-algebra $ {\mathcal A} $ and abelian $ n $-tuples 
    $ \underline{x}=(x_1,\dots,x_n) $ and $ \underline{y}=(y_1,\dots,y_n) $ of positive elements in $ {\mathcal A}. $
    
    In section \ref{section: expectation values} we consider functionals which are given as the
    expectation value of a unit vector in a Hilbert space,
    and derive what we term Jensen's inequality for expectation values. 
    
    In section \ref{last section} we finally prove some related results for matrix functions of several variables, but
    with respect to the weak majorization for matrices.

    \section{Inequalities under a positive functional}
    
    Let $ {\mathcal C} $ be a separable abelian $ C^* $-subalgebra of a $ C^* $-algebra $ {\mathcal A}, $ and let
    $ \varphi $ be a positive functional on $ {\mathcal A} $ such that $ {\mathcal C} $ is contained in the
    centralizer $ {\mathcal A}^\varphi. $ The subalgebra is of the form $ {\mathcal C}=C_0(S) $
    for some locally compact metric space $ S, $ and by the Riesz representation theorem there is a finite
    Radon measure $ \mu_\varphi $ on $ S $ such that
    \[
    \varphi(y)=\int_S y(s)\,d\mu_\varphi(s)\qquad y\in{\mathcal C}=C_0(S).
    \]    
    For each positive element $ x $ in the multiplier algebra
    $ M(\mathcal A) $ we have
    \[
    0\le\varphi(xy)=\varphi(y^{1/2}xy^{1/2})\le\|x\|\varphi(y)\qquad y\in{\mathcal C}_+.
    \]
    The functional $ y\to\varphi(xy) $ on $ {\mathcal C} $ consequently defines a Radon measure on $ S $ which is
    dominated by a multiple of $ \mu_\varphi, $ and it is therefore given by a unique element $ \Phi(x) $ in
    $ L^\infty(S,\mu_\varphi). $ By linearization this defines a conditional expectation\footnote{This
    is a slight abuse of language since the range is not a subalgebra of $ M({\mathcal A}), $ but $ \Phi $
    is positive, linear and $ \Phi(xy)=\Phi(x)y $ for $ x\in M({\mathcal A}) $ and $ y\in {\mathcal C}.$}
    of the multiplier algebra
    \begin{gather}\label{conditional expectation}
    \Phi\colon M({\mathcal A})\to L^\infty(S,\mu_\varphi)
    \end{gather}
    such that
    \[
    \int_S z(s)\Phi(x)(s)\, d\mu_\varphi(s)=\varphi(zx)\qquad z\in {\mathcal C},\, x\in M({\mathcal A}).
    \]
    In particular, $ \Phi(z)(s)=z(s) $ almost everywhere in $ S $ for each $ z\in {\mathcal C}, $
    cf. \cite{kn:lieb:2002, kn:hansen:2003:2, kn:hansen:2003:3}. 
    
    The following result is, although not explicitly stated, essentially proved in \cite{kn:hansen:2003:3}
    and follows by inspection of the proof of Theorem 4.1 in the reference.
    
    \begin{theorem}    
    Let $ f\colon\underline{I}\to{\mathbf R} $ be a continuous function defined on a cube 
    $ \underline{I}=I_1\times\cdots\times I_n $ and let $ \underline{x}=(x_1,\dots,x_n)$ be an
    abelian $ n $-tuple in $ {\mathcal A} $ contained in the domain of $ f. $
    If $ f $ is concave then
    \[
    \Phi(f(\underline{x}))\le f(\Phi(x_1),\dots,\Phi(x_n))
    \]
    almost everywhere,
    where $ \Phi\colon M({\mathcal A})\to L^\infty(S,\mu_\varphi) $ is the conditional expectation 
    in (\ref{conditional expectation}).    
    \end{theorem}
    
    If in addition $ y=(y_1,\dots,y_n) $ is an $ n $-tuple in $ {\mathcal C} $ contained in the domain of $ f, $ 
    and $ f $ is also separately increasing, then $ \underline{x}\le\underline{y} $ implies
    \begin{gather}\label{monotonicity under Phi}
    \Phi(f(\underline{x}))\le f(\Phi(x_1),\dots,\Phi(x_n))\le f(\Phi(y_1),\dots,\Phi(y_n))=f(\underline{y})
    \end{gather}
    almost everywhere, and consequently
    \[
    \varphi(f(\underline{x}))=\varphi(\Phi(f(\underline{x}))\le\varphi(f(\underline{y})).
    \]
    The monotonicity under $ \varphi $ is thus a trivial consequence of the assertion proved
    in equation (\ref{monotonicity under Phi}), and it therefore comes as no surprise that we can find examples where
    (\ref{monotonicity under Phi}) is violated, but we still have monotonicity under $ \varphi. $
    
    \begin{example}
    Let $ {\mathcal A}=M_2 $ be the $ C^* $-algebra of $ 2\times 2 $ matrices, and let $ {\mathcal C} $ be the abelian
    subalgebra of the diagonal matrices. We consider the matrices
    \[
    x=\begin{pmatrix}
    c & c\\
    c & c
    \end{pmatrix}\quad\mbox{and}\quad
    y=\begin{pmatrix}
    t & 0\\
    0 & \lambda t
    \end{pmatrix}\qquad 0<c<t
    \]
    and notice that $ 0\le x< y $ for $ \lambda> c(t-c)^{-1}. $ We choose the ordinary
    trace as the positive functional and obtain
    \[
    \Phi(x^2)=\begin{pmatrix}
    2c^2 & 0\\
    0 & 2c^2
    \end{pmatrix}\not\leq\,\,
    \begin{pmatrix}
    t^2 & 0\\
    0   & \lambda^2 t^2
    \end{pmatrix}=y^2\qquad\mbox{for}\quad t< c\sqrt{2},
    \]
    while as expected $ \tr x^2=4 c^2<t^2(1+c^2(t-c)^{-2})\le t^2(1+\lambda^2)=\tr y^2. $     
    \end{example}
    
    The theory of operator means \cite{kn:kubo:1980} provides us with a key tool to obtain the implication in   
    (\ref{monotone trace function}) for functions which are neither convex nor concave,
    provided we assume that all of the elements $ x_1,\dots,x_n $ and $ y_1,\dots,y_n $
    are contained in the centralizer $ {\mathcal A}^\varphi. $
    
    \begin{theorem}
    Let $ \varphi $ be a positive functional on a $ C^* $-algebra $ {\mathcal A}, $ and let
    $ \underline{x}=(x_1,\dots,x_n) $ and $ \underline{y}=(y_1,\dots,y_n) $ be abelian $ n $-tuples 
    of elements in the centralizer
    $ {\mathcal A}^\varphi. $ If $ 0\le x_i\le y_i $ for $ i=1,\dots,n $ then
    \[
    \varphi(x_1^{p_1}\cdots x_n^{p_n})\le\varphi(y_1^{p_1}\cdots y_n^{p_n})
    \]
    for arbitrary non-negative exponents $ p_1,\dots,p_n. $ 
    \end{theorem}
    
    \proof The geometric mean $ x\# y $ is defined for positive invertible elements $ x,y\in {\mathcal A} $ 
    by setting
    \[
    x\# y=x^{1/2}(x^{-1/2}yx^{-1/2})^{1/2}x^{1/2}
    \]
    and since 
    \[
    x^{1/2}(x^{-1/2}yx^{-1/2})^{1/2}x^{1/2}=   
    \frac{1}{2\pi}\int_0^\infty 2(x^{-1}+\lambda y^{-1})^{-1}\lambda^{-1/2}\,d\lambda,
    \]
    it follows that $ x\# y $ is increasing in each variable. In fact,    
    it can be extended to positive elements $ x,y\ge 0 $ in $ {\mathcal A} $ and becomes a concave and separately 
    increasing function of the pair
    $ (x,y), $ cf. \cite{kn:kubo:1980}.  Since $ x_1 $ and $ x_2 $ commute we therefore obtain
    \[
    x_1^{1/2}x_2^{1/2}=x_1\# x_2 \le y_1\# y_2=y_1^{1/2}y_2^{1/2}.
    \]
    We first note that $ x_3^{1/2}\le y_3^{1/2} $ by the Löwner-Heinz inequality. Furthermore, $ x_3^{1/2} $
     commutes with
    $ x_1^{1/2}x_2^{1/2} $ (with the same statement for the $ y $'s). 
    We may therefore apply the above procedure once more to the
    abelian pairs $ (x_1^{1/2}x_2^{1/2}, x_3^{1/2}) $ and $ (y_1^{1/2}y_2^{1/2}, y_3^{1/2}) $ to get
    \[
    x_1^{1/4}x_2^{1/4}x_3^{1/4}\le y_1^{1/4}y_2^{1/4}y_3^{1/4}
    \]
    and by induction we finally obtain
    \[
    x_1^{1/2^{n-1}}\cdots x_n^{1/2^{n-1}}\le y_1^{1/2^{n-1}}\cdots y_1^{1/2^{n-1}}.
    \]
    For any continuous and increasing function $ f $ defined on the positive half-axis we
    have
    \[
    0\le x\le y\quad\Rightarrow\quad \varphi(f(x))\le\varphi(f(y))
    \]
    for elements $ x $ and $ y $ in the centralizer $ {\mathcal A}^\varphi, 
    $ cf. \cite[Theorem 4.2]{kn:hansen:2003:3}
    (note that the compatibility condition in the reference is void for functions of one variable); cf. also
    \cite{kn:hansen:1995, kn:brown:1990, kn:bernstein:1988}. Setting $ f(t)=t^{2^{n-1}\cdot N} $ we therefore
    obtain
    \[
    \varphi(x_1^N\cdots x_n^N)\le\varphi(y_1^N\cdots y_n^N)
    \]
    for arbitrary $ N>0. $ Possibly by first applying Löwner-Heinz inequality for each entry we thus have
    \[
    \varphi(x_1^{\alpha_1 N}\cdots x_n^{\alpha_n N})\le\varphi(y_1^{\alpha_1 N}\cdots y_n^{\alpha_n N})
    \]
    for $ \alpha_1,\dots,\alpha_n\in[0,1] $ and $ N>0, $ and since any set of positive
    exponents $ (p_1,\dots,p_n) $ can be written in this form the assertion follows.
    \qed
    
    \section{Jensen's inequality for expectation values}\label{section: expectation values}
    
    Recall that a continuous field $ t\to a_t $ of operators on a Hilbert space $ H $ defined on a locally
    compact Hausdorff space $ T $ equipped with a Radon measure $ \nu $ is said to be a unital column field if
    \[
    \int_T a_t^*a_t\, d\nu(t)=1,
    \]
    cf. \cite{kn:hansen:2003:2}.
    
    \begin{theorem}\label{Jensen's inequality for expectation values}
    Let $ f:\underline{I}\to {\mathbf R} $ be a continuous convex function of $ n $ variables defined on a cube,
    and let $ t\to a_t\in B(H) $ be a unital column field on a locally compact Hausdorff space $ T $ with a Radon
    measure $ \nu. $ If $ t\to\underline{x}_t $ is a bounded continuous field on $ T $ of abelian 
    $ n $-tuples of operators on $ H $ in the domain of $ f, $ then
    \begin{gather}\label{vector inequality}
    f\bigl((y_1\xi\mid\xi),\dots,(y_n\xi\mid\xi)\bigr)\le\left( 
    \int_T a_t^* f(\underline{x}_t) a_t\, d\nu(t)\xi\mid\xi\right) 
    \end{gather}
    for any unit vector $ \xi\in H, $ where the $ n $-tuple $ \underline{y} $ is defined by setting
    \[
    \underline{y}=(y_1,\dots,y_n)=\int_T a_t^* \underline{x}_t a_t\,d\nu(t).
    \]
    \end{theorem}
    
    \begin{proof}
    Set $ \underline{x}_t=(x_{1t},\dots,x_{nt}) $ for $ t\in T $ and consider the spectral resolutions
    \[
    x_{it}=\int\lambda\, dE_{it}(\lambda)\qquad i=1,\dots,n;\, t\in T.
    \] 
    There is to each unit vector $ \xi\in H $ a positive measure $ \mu_\xi $ on $ \underline{I} $
    such that
    \[
    \mu_\xi(S_1\times\cdots\times S_n)=\int_T (E_{1t}(S_1)\cdots E_{nt}(S_n) a_t\xi\mid a_t\xi)\, d\nu(t)
    \]
    for Borel sets $ S_1\subseteq I_1,\dots,S_n\subseteq I_n $ and since the column field $ t\to a_t $ 
    is unital, we obtain
    that $ \mu_\xi $ is a probability measure. It satisfies
    \[\label{representation of operator functions}
    \int_T \bigl(g(\underline{x}_t) a_t\xi\mid a_t\xi\bigr)\, d\nu(t)
    =\int_{\underline{I}} g(\underline{s})\,d\mu_\xi(\underline{s})\qquad\underline{s}=(s_1,\dots,s_n)
    \]
    for any continuous function $ g:\underline{I}\to{\mathbf R}. $ In particular 
    (putting $ g_i(\underline{s})=s_i) $ we obtain
    \[
    \int_T \bigl(x_{it} a_t\xi\mid a_t\xi\bigr)\, d\nu(t)=\int_{\underline{I}} s_i\,d\mu_\xi(\underline{s})
    \qquad i=1,\dots,n.
    \]
    We thus obtain
    \[
    \begin{array}{rl}
    &f\bigl((y_1\xi\mid\xi),\dots, (y_n\xi\mid\xi)\bigr)\\[2ex]
    =&\displaystyle f\left(\left(\int_T a_t^* x_{1t} a_t \,d\nu(t)\xi\mid\xi\right),
    \dots, \left(\int_T a_t^* x_{nt} a_t \,d\nu(t)\xi\mid\xi\right)\right)\\[3ex]
    =&\displaystyle f\left(\int_T (x_{1t} a_t\xi\mid a_t\xi)\,d\nu(t),
    \dots,\int_T (x_{nt} a_t\xi\mid a_t\xi)\,d\nu(t)\right)\\[3ex]
    =&\displaystyle f\left(\int_{\underline{I}} s_1\,d\mu_\xi(\underline{s}),
    \dots,\int_{\underline{I}} s_n\,d\mu_\xi(\underline{s})\right)\\[3ex]
    \le&\displaystyle\int_{\underline{I}} f(s_1,\dots,s_n)\, d\mu_\xi(\underline{s})
    =\displaystyle\int_T \bigl(f(\underline{x}_t) a_t \xi\mid a_t \xi\bigr)\, d\nu(t)\\[3ex]  
    =&\displaystyle\left(\int_T a_t^* f(\underline{x}_t) a_t\,d\nu(t) \xi\mid\xi\right),   
    \end{array}
    \]
    where we used spectral theory and the convexity of $ f. $ 
    \end{proof} 
    
     \begin{remark}\label{remark}
    If we choose $ \nu $ as a probability measure on $ T, $ then the trivial field $ a_t=1 $
    for $ t\in T $ is unital and (\ref{vector inequality}) takes the form
    \[
    f\left(\left(\int_T x_{1t} \,d\nu(t)\xi\mid\xi\right),
    \dots, \left(\int_T x_{nt} \,d\nu(t)\xi\mid\xi\right)\right)\le
    \left(\int_T f(\underline{x}_t)\,d\nu(t) \xi\mid\xi\right)
    \]
    for fields of abelian $ n $-tuples $ \underline{x}_t=(x_{1t},\dots,x_{nt}) $ and unit vectors $ \xi. $
    By choosing $ \nu $ as an atomic measure with one atom we get a version 
    \begin{gather}\label{mond and pecaric}
    f\bigl(\left(x_1\xi\mid\xi\right),\dots, \left(x_n\xi\mid\xi\right)\bigr)\le
    \bigl(f(\underline{x})\xi\mid\xi\bigr)
    \end{gather}
    of the Jensen inequality by Mond and Pe{\v{c}}ari{\'{c}} \cite{kn:mond:1993:1}.
    \end{remark}
    
    One may generalize Theorem \ref{Jensen's inequality for expectation values} 
    by using the notions and methods developed in the proof of \cite[Theorem 3.1]{kn:hansen:2003:3}
    to obtain the following result.
    
    \begin{theorem}
     Let $ {\mathcal C} $ be a separable abelian $ C^* $-subalgebra of a $ C^* $-algebra $ {\mathcal A}, $ let
    $ \varphi $ be a positive functional on $ {\mathcal A} $ such that $ {\mathcal C} $ is contained in the
    centralizer $ {\mathcal A}^\varphi $ and let 
    \[
    \Phi\colon M({\mathcal A})\to L^\infty(S,\mu_\varphi) 
    \]
    be the conditional expectation defined in (\ref{conditional expectation}). Let furthermore
    $ f:\underline{I}\to {\mathbf R} $ be a continuous convex function of $ n $ variables defined on a cube,
    and let $ t\to a_t\in M({\mathcal A}) $ be a unital column field on a locally compact Hausdorff space $ T $
    with a Radon
    measure $ \nu. $ If $ t\to\underline{x}_t $ is a bounded, weak* measurable field on $ T $ of abelian 
    $ n $-tuples in $ {\mathcal A} $ in the domain of $ f, $ then
    \[
    f(\Phi(y_1),\dots,\Phi(y_n))\le \Phi\left(\int_T a_t^*f(\underline{x}_t) a_t\,d\nu(t)\right)
    \]
    almost everywhere,
    where the $ n $-tuple $ \underline{y} $ in $ M({\mathcal A}) $ is defined by setting
    \[
    \underline{y}=(y_1,\dots,y_n)=\int_T a_t^* \underline{x}_t a_t\,d\nu(t).    
    \]
    \end{theorem}
    
    Note that Theorem \ref{Jensen's inequality for expectation values} follows from the preceding theorem
    by choosing $ \varphi $ as the trace and letting $ {\mathcal C} $ be the $ C^* $-algebra generated by the
    orthogonal projection on the vector $ \xi. $

     \section{Weak majorization for matrices}\label{last section}   

    We consider a Hilbert space $ H $ of finite dimension $ m $ and introduce for any self-adjoint operator 
    $ x\in B(H) $ the $ m $-tuple $ (x_{[1]},\dots,x_{[m]}) $ of eigenvalues of $ x $ counted with
    multiplicity and ordered in a decreasing sequence. The notion of weak majorization for matrices was considered
    by Ando \cite{kn:ando:1994} and Bhatia \cite{kn:bhatia:1997}. 
    
    \begin{definition} Let $ x $ and $ y $ be self-adjoint operators on a Hilbert space $ H $ of finite dimension $ m. $
    We say that $ x $ is weakly majorized by $ y, $ and we write $ x\prec_{\mbox{\tiny w}} y $ if
    \[
    \sum_{i=1}^k x_{[i]}\le\sum_{i=1}^k y_{[i]}
    \]
    for $ k=1,\dots,m. $
    \end{definition}
    
    The following result is known as Ky Fan's maximum principle, cf. Bhatia \cite[p. 35]{kn:bhatia:1997}.
    
    \begin{lemma}\label{lemma: majorization}
    Let $ x $ be a self-adjoint operator on a Hilbert space $ H $ of finite dimension $ m, $ and take 
    a natural number $ k\le m. $ Then
    \[
    \sum_{i=1}^k (xu_i\mid u_i)\le\sum_{i=1}^k x_{[i]}
    \]
    for any orthonormal set $ (u_1,\dots,u_k) $ of vectors in $ H. $
    \end{lemma}
    
   \begin{theorem} Let $ f:\underline{I}\to {\mathbf R} $ be a convex function of $ n $ variables 
   defined on a cube, let $ H $ be a Hilbert space of finite dimension
    and let $ t\to a_t\in B(H) $ be a unital column field on a locally compact Hausdorff space $ T $
    with a Radon measure $ \nu. $ If $ t\to\underline{x}_t $ is a bounded continuous field on $ T $ of abelian 
    $ n $-tuples of operators on $ H $ in the domain of $ f, $ then the inequality    
    \[
    f\left(\int_T a_t^* \underline{x}_t a_t\,d\nu(t)\right) \prec_{\mbox{\tiny w}}
    \int_T a_t^* f(\underline{x}_t) a_t\, d\nu(t)
    \]
    is valid provided 
    \[
    \underline{y}=(y_1,\dots,y_n)=\int_T a_t^* \underline{x}_t a_t\,d\nu(t)
    \]
    is an abelian $ n $-tuple.    
    \end{theorem}
    
    Note that the assumption of $ \underline{y} $ being an abelian $ n $-tuple
    is void for functions of one variable.\vskip 2ex
    
    \begin{proof}    
    For $ m=\dim H $ we choose an orthonormal $ m $-tuple $ (u_1,\dots,u_m) $ of 
    common eigenvectors for the commuting matrices 
    $ y_1,\dots,y_n $ in such a way that the corresponding eigenvalues of $ f(\underline{y}) $ are ordered in a 
    decreasing sequence. We then obtain
    \[
    \begin{array}{rl}
    &\displaystyle\sum_{i=1}^k f\left(\int_T a_t^* \underline{x}_t a_t\,d\nu(t)\right)_{[i]}
    =\sum_{i=1}^k f(\underline{y})_{[i]}
    =\sum_{i=1}^k \bigl(f(\underline{y})u_i\mid u_i\bigr)\\[2ex]
    =&\displaystyle\sum_{i=1}^k 
    f\bigl((y_1u_i\mid u_i),\dots, (y_n u_i\mid u_i)\bigr)\\[2ex]
    \le&\displaystyle\sum_{i=1}^k\left(\int_T a_t^* f(\underline{x}_t) a_t\,d\nu(t) u_i\mid u_i\right)\\[2ex]  
    \le&\displaystyle\sum_{i=1}^k 
    \left(\int_T a_t^* f(\underline{x}_t) a_t\,d\nu(t)\right)_{[i]}\qquad k=1,\dots, m,    
    \end{array}
    \]
    where we used spectral theory, Jensen's inequality for expectation values (\ref{vector inequality})
    and Lemma \ref{lemma: majorization}.
    \end{proof}
    
    \begin{remark}  If we choose a probability measure $ \nu $ and the trivial field $ a_t=1 $
    we obtain the inequality
    \[
    f\left(\int_T \underline{x}_t \,d\nu(t)\right) \prec_{\mbox{\tiny w}}
    \int_T f(\underline{x}_t) \, d\nu(t),
    \]
    provided the integral $ \displaystyle\int_T \underline{x}_t \,d\nu(t) $ is an abelian $ n $-tuple.        
    \end{remark}
    
    The following result is a generalization to functions of several variables of a theorem by 
    Aujla and Silva \cite[Theorem 2.3]{kn:aujla:2003} for functions of one variable.
    
    \begin{corollary}
    Let $ f:\underline{I}\to {\mathbf R} $ be a convex function of $ n $ variables defined on a cube
    $ \underline{I}, $ and let $ H $ be a Hilbert space of finite dimension. Then
    \[
    f(\lambda \underline{x}+(1-\lambda)\underline{y})\prec_{\mbox{\tiny w}}
    \lambda f(\underline{x})+(1-\lambda)f(\underline{y})\qquad\lambda\in[0,1]
    \]
    for compatible $ n $-tuples $ \underline{x} $ and $ \underline{y} $ of operators on $ H $
    in the domain of $ f. $
    \end{corollary}
    
    \begin{proof} Since $ \underline{x} $ and $ \underline{y} $ are compatible, the line through 
    $ \underline{x} $ and $ \underline{y} $ consists of abelian $ n $-tuples. The result therefore
    follows from the preceding remark by choosing a suitable atomic probability measure $ \nu. $ 
    \end{proof}
    
    \begin{theorem}
    Let $ f:\underline{I}\to {\mathbf R} $ be a convex and separately increasing
    function of $ n $ variables, and let $ H $ be a Hilbert space of finite dimension. Then 
    \[
    \underline{x}\le\underline{y}\quad\Rightarrow\quad
    f(\underline{x})\prec_{\mbox{\tiny w}} f(\underline{y})
    \]
    for abelian $ n $-tuples $ \underline{x}=(x_1,\dots,x_n) $ and $ \underline{y}=(y_1,\dots,y_n) $
    of operators on $ H $ in the domain of $ f. $
    \end{theorem}
    
    \begin{proof}
    For $ m=\dim H $ we choose an orthonormal $ m $-tuple $ (u_1,\dots,u_m) $ 
    of common eigenvectors for the commuting
    matrices $ x_1,\dots,x_n $ in such a way that the corresponding eigenvalues for $ f(\underline{x}) $
    are ordered in a decreasing sequence. We then obtain
    \[
    \begin{array}{l}\displaystyle
    \sum_{i=1}^k f(\underline{x})_{[i]}
    =\sum_{i=1}^k \bigl(f(\underline{x})u_i\mid u_i\bigr)
    =\sum_{i=1}^k f\bigl((x_1u_i\mid u_i),\dots, (x_nu_i\mid u_i)\bigr)\\[2.5ex]
    \le \displaystyle\sum_{i=1}^k 
    f\bigl((y_1u_i\mid u_i),\dots, (y_nu_i\mid u_i)\bigr)
    \le\sum_{i=1}^k \bigl(f(\underline{y})u_i\mid u_i\bigr)\\[2.5ex]
    \le\displaystyle\sum_{i=1}^k 
    f(\underline{y})_{[i]}\qquad k=1,\dots,m,    
    \end{array}
    \]
    where we used spectral theory, the monotonicity and convexity of $ f, $ the inequality
    in (\ref{mond and pecaric}) and Lemma \ref{lemma: majorization} respectively.
    \end{proof}
    
    \nocite{kn:pedersen:2003}
    \nocite{kn:marshall:1979}
    
     \bibliographystyle{plain}

      \vfill

      {\small\noindent Frank Hansen: Institute of Economics, University
       of Copenhagen, Studiestraede 6, DK-1455 Copenhagen K, Denmark.}

      \end{document}